\documentclass[11pt]{article}
\usepackage{amsbsy,amsfonts,amsmath,amssymb,amsthm,bbm,enumerate}

\numberwithin{equation}{section}

\theoremstyle{plain}
\newtheorem{theorem}{Theorem}

\newtheorem{lemma}[theorem]{Lemma}

\newcommand{\bb}{\underline{b}}

\newcommand{\cC}{{\cal C}}

\newcommand{\db}{\rightrightarrows}

\newcommand{\dpst}{\displaystyle}

\newcommand{\lbeq}[1]{\label{eq:#1}}

\newcommand{\mN}{{\mathbb N}}
\newcommand{\mP}{{\mathbb P}}

\newcommand{\mZ}{{\mathbb Z}}
\newcommand{\nn}{\nonumber}

\newcommand{\piv}{{\tt piv}}
\newcommand{\Proof}[1]{\smallskip\noindent\emph{#1.}~ }
\newcommand{\QED}{\hspace*{\fill}\rule{7pt}{7pt}\smallskip}

\newcommand{\refeq}[1]{(\ref{eq:#1})}
\newcommand{\sss}{\scriptscriptstyle}

\newcommand{\tb}{\overline{b}}
\newcommand{\tbp}{\tb^{\raisebox{-2.4pt}{\scriptsize$\prime$}}}
\newcommand{\thetac}{\theta^\times}
\newcommand{\thetap}{\theta^{\sss|\,|}}

\newcommand{\Zd}{\mZ^d}
\newcommand{\xic}{\xi^\times}
\newcommand{\xip}{\xi^{\sss|\,|}}
\newcommand{\Zp}{\mZ_+}

\newcommand{\ovec}  {\boldsymbol{o}}

\newcommand{\uvec}  {\boldsymbol{u}}
\newcommand{\vvec}  {\boldsymbol{v}}
\newcommand{\wvec}  {\boldsymbol{w}}
\newcommand{\xvec}  {\boldsymbol{x}}
\newcommand{\yvec}  {\boldsymbol{y}}

\title{\textbf{Unpublished manuscript}\\[1pc]
Diagrammatic bounds on the lace-expansion coefficients\\
for oriented percolation}
\author{Akira~Sakai\footnote{Department of Mathematical Sciences,
University of Bath, UK.}}
\date{July 25, 2007}

\begin{document}
\maketitle

In this note, we provide a complete proof of \cite[Proposition~3.3]{cs?}.  For
notational convenience, we use bold letters to denote vertices in $\mZ^{d+1}$,
e.g., $\ovec\equiv(o,0)$ and $\xvec$; if necessary, we denote the spatial and
temporal components of a given vertex $\vvec$ by $\sigma_{\vvec}$ and
$\tau_{\vvec}$ respectively: $\vvec=(\sigma_{\vvec},\tau_{\vvec})$.  To
identify the starting and terminal points, we write, e.g.,
$\varphi_p(\vvec;\xvec)=\mP_p(\vvec\to\xvec)$ and abbreviate it to
$\varphi_p(\xvec)$ if $\vvec=\ovec$; in particular,
$\varphi_p(\vvec;\xvec)=\varphi_p(\xvec-\vvec)$ if the model is
translation-invariant.  Let $\piv(\vvec,\xvec)$ denote the (random) set of
pivotal bonds for $\{\vvec\to\xvec\}$.

\section{Bounds in terms of two-point functions}
In this section, we prove bounds on $\pi_p^{\sss(N)}(\xvec)$ and
$\Pi_p^{\sss(N)}(\xvec)$, for fixed $\xvec$, in terms of two-point functions.
To prove these bounds, we do not have to assume translation-invariance.

Recall that the lace-expansion coefficients $\pi_p^{\sss(N)}(\xvec)$ and
$\Pi_p^{\sss(N)}(\xvec)$ for $N\ge1$ are defined in terms of the event
\begin{align}
\tilde E_{\vec b_N}^{\sss(N)}(\xvec)=\{\ovec\db\bb_1\}\cap\bigcap_{i=1}^N
 E\big(b_i,\bb_{i+1};\tilde\cC^{b_i}(\tb_{i-1})\big),
\end{align}
where $\vec b_N=(b_1,\dots,b_N)$ is an ordered set of bonds and
\begin{align}
E(b,\xvec;\cC)&=\{b\to\xvec\in\cC\}\cap\big\{\nexists\,b'\in\piv(\tb,\xvec)
 \text{ satisfying }\bb'\in\cC\big\},\\
\tilde\cC^b(\vvec)&=\{\xvec\in\mZ^{d+1}:\vvec\to\xvec\text{ without using }b\}.
\end{align}

\begin{lemma}\label{lem:piN-BKappl}
\begin{align}\lbeq{pi0-BKappl}
\pi_p^{\sss(0)}(\xvec)\equiv\mP_p(\ovec\db\xvec)\le\delta_{\xvec,\ovec}+(q_p*
 \varphi_p)(\xvec)^2,
\end{align}
and, for $N\ge1$,
\begin{align}\lbeq{piN-BKappl}
\pi_p^{\sss(N)}(\xvec)\equiv\sum_{\vec b_N}\mP_p\big(\tilde E_{\vec b_N}^{\sss
 (N)}(\xvec)\big)\le\sum_{\substack{\uvec_1,\dots,\uvec_{N+1}\\ \vvec_1,\dots,
 \vvec_{N+1}\\ (\uvec_{N+1}=\vvec_{N+1}=\xvec)}}\varphi_p(\uvec_1)\,\varphi_p
 (\uvec_1;\vvec_1)\,\varphi_p(\vvec_1)\prod_{i=1}^N\Xi_p(\uvec_i,\vvec_i;
 \uvec_{i+1},\vvec_{i+1}),
\end{align}
where
\begin{gather}
\Xi_p(\uvec,\vvec;\uvec',\vvec')=\big(\xip_p(\uvec,\vvec;\uvec',\vvec')+\xic_p
 (\uvec,\vvec;\uvec',\vvec')\big)\varphi_p(\uvec';\vvec')/2^{\delta_{\uvec',
 \vvec'}},\lbeq{xi-def}\\
\begin{cases}
\xip_p(\uvec,\vvec;\uvec',\vvec')=(q_p*\varphi_p)(\uvec;\uvec')\,(q_p*
 \varphi_p)(\vvec;\vvec'),\\
\xic_p(\uvec,\vvec;\uvec',\vvec')=(q_p*\varphi_p)(\uvec;\vvec')\,(q_p*
 \varphi_p)(\vvec;\uvec').
\end{cases}
\end{gather}
\end{lemma}

\Proof{Proof} Since \refeq{pi0-BKappl} is already proved in \cite[(3.18)]{cs?},
it remains to show \refeq{piN-BKappl}.  By definition, we can easily see that
\begin{align}\lbeq{Esup}
E(b,\xvec;\tilde\cC^b(\yvec))\subset\{\yvec\to\xvec\}\circ\{b\to\xvec\},
\end{align}
where $E_1\circ E_2$ is the event that $E_1$ and $E_2$ occur bond-disjointly
(i.e., $E_1$ occurs on some bond set $B$ and $E_2$ occurs on $B^\text{c}$).
Similarly,
\begin{align}
\{\ovec\db\vvec\}\cap\{\ovec\to\xvec\}&\subset\bigcup_{\uvec}\big\{\{
 \ovec\to\uvec\to\vvec\}\circ\{\ovec\to\vvec\}\circ\{\uvec\to\xvec\}
 \big\},\lbeq{dbellsup}\\
E(b,\vvec;\tilde\cC^b(\yvec))\cap\{\tb\to\xvec\}&\subset\bigcup_{\uvec:
 \tau_{\uvec}>\tau_{\bb}}\Big\{\big\{\{\yvec\to\uvec\to\vvec\}\circ\{b
 \to\vvec\}\circ\{\uvec\to\xvec\}\big\}\nn\\
&\hspace{3.5pc}\cup\big\{\{\yvec\to\vvec\}\circ\{b\to\uvec\to\vvec\}\circ\{
 \uvec\to\xvec\}\big\}\Big\}.\lbeq{Eellsup}
\end{align}
To prove \refeq{piN-BKappl}, we use \refeq{Esup}--\refeq{Eellsup} and the BK
inequality and pay attention to which event depends on which time interval.
For example, by \refeq{Esup},
\begin{align}
\tilde E_{\vec b_N}^{\sss(N)}(\xvec)\subset\tilde E_{\vec b_{N-1}}^{\sss(N-1)}
 (\bb_N)\cap\{\tb_{N-1}\to\xvec\}\circ\{b_N\to\xvec\}.
\end{align}
Since $\tilde E_{\vec b_{N-1}}^{\sss(N-1)}(\bb_N)$ depends only on bonds before
time $\tau_{\bb_N}$, we can use the BK inequality to obtain
\begin{align}\lbeq{piN-1stbd}
\sum_{b_N}\mP_p\big(\tilde E_{\vec b_N}^{\sss(N)}(\xvec)\big)\le\sum_{\vvec_N}
 \mP_p\Big(\tilde E_{\vec b_{N-1}}^{\sss(N-1)}(\vvec_N)\cap\{\tb_{N-1}\to\xvec
 \}\Big)\,(q_p*\varphi_p)(\vvec_N;\xvec).
\end{align}
Then, by \refeq{Eellsup} and the BK inequality and using the Markov property,
we obtain
\begin{align}\lbeq{piN-2ndbd}
&\sum_{b_{N-1}}\mP_p\Big(\tilde E_{\vec b_{N-1}}^{\sss(N-1)}(\vvec_N)\cap\{
 \tb_{N-1}\to\xvec\}\Big)\nn\\
&\le\sum_{\substack{\vvec_{N-1},\uvec_N\\ (\tau_{\vvec_{N-1}}<\tau_{\uvec_N})}}
 \bigg(\mP_p\Big(\tilde E_{\vec b_{N-2}}^{\sss(N-2)}(\vvec_{N-1})\cap\{\tb_{N
 -2}\to\uvec_N\}\Big)(q_p*\varphi_p)(\vvec_{N-1};\vvec_N)\\
&\hspace{6pc}+\mP_p\Big(\tilde E_{\vec b_{N-2}}^{\sss(N-2)}(\vvec_{N-1})\cap
 \{\tb_{N-2}\to\vvec_N\}\Big)(q_p*\varphi_p)(\vvec_{N-1};\uvec_N)\bigg)
 \varphi_p(\uvec_N;\vvec_N)\,\varphi_p(\uvec_N;\xvec).\nn
\end{align}
Since $\tau_{\uvec_N}\le\tau_{\vvec_N}<\tau_{\xvec}$ (due to
$(q_p*\varphi_p)(\vvec_N;\xvec)$ in \refeq{piN-1stbd} and
$\varphi_p(\uvec_N;\vvec_N)$ in \refeq{piN-2ndbd}), we can replace the last
term in \refeq{piN-2ndbd} by $(q_p*\varphi_p)(\uvec_N;\xvec)$, using the
trivial inequality
\begin{align}\lbeq{varphi-trivineq}
\varphi_p(\uvec;\xvec)\le(q_p*\varphi_p)(\uvec;\xvec)\qquad(\uvec\ne\xvec).
\end{align}
Summarizing these bounds, we have
\begin{align}\lbeq{piN-3rdbd}
\sum_{b_{N-1},b_N}\mP_p\big(\tilde E_{\vec b_N}^{\sss(N)}(\xvec)\big)\le\sum_{
 \substack{\vvec_{N-1},\uvec_N,\vvec_N\\ (\tau_{\vvec_{N-1}}<\tau_{\uvec_N})}}
 &\bigg(\mP_p\Big(\tilde E_{\vec b_{N-2}}^{\sss(N-2)}(\vvec_{N-1})\cap\{
 \tb_{N-2}\to\uvec_N\}\Big)(q_p*\varphi_p)(\vvec_{N-1};\vvec_N)\nn\\
&+\mP_p\Big(\tilde E_{\vec b_{N-2}}^{\sss(N-2)}(\vvec_{N-1})\cap\{\tb_{N-2}
 \to\vvec_N\}\Big)(q_p*\varphi_p)(\vvec_{N-1};\uvec_N)\bigg)\nn\\
&\hspace{-1pc}\times\varphi_p(\uvec_N;\vvec_N)\,\Xi_p(\uvec_N,\vvec_N;
 \xvec,\xvec).
\end{align}
Using \refeq{piN-2ndbd}--\refeq{varphi-trivineq} again, but with different
variables, we obtain
\begin{align}\lbeq{piN-4thbd}
\sum_{b_{N-2},b_{N-1},b_N}\mP_p\big(\tilde E_{\vec b_N}^{\sss(N)}(\xvec)\big)
 &\le\sum_{\substack{\uvec_{N-1},\uvec_N\\ \vvec_{N-2},\vvec_{N-1},\vvec_N\\
 (\tau_{\vvec_{N-2}}<\tau_{\uvec_{N-1}})}}\bigg(\mP_p\Big(\tilde E_{\vec b_{
 N-3}}^{\sss(N-3)}(\vvec_{N-2})\cap\{\tb_{N-3}\to\uvec_{N-1}\}\Big)(q_p*
 \varphi_p)(\vvec_{N-2};\vvec_{N-1})\nn\\
&\qquad\qquad\qquad+\mP_p\Big(\tilde E_{\vec b_{N-3}}^{\sss(N-3)}(\vvec_{N-2})
 \cap\{\tb_{N-3}\to\vvec_{N-1}\}\Big)(q_p*\varphi_p)(\vvec_{N-2};\uvec_{N-1})
 \bigg)\nn\\
&\qquad\qquad\times\varphi_p(\uvec_{N-1};\vvec_{N-1})\,\Xi_p(\uvec_{N-1},
 \vvec_{N-1};\uvec_N,\vvec_N)\,\Xi_p(\uvec_N,\vvec_N;\xvec,\xvec).
\end{align}
We repeat this procedure until we arrive at
\begin{align}\lbeq{piN-5thbd}
\sum_{\vec b_N}\mP_p\big(\tilde E_{\vec b_N}^{\sss(N)}(\xvec)\big)\le\sum_{
 \substack{\uvec_2,\dots,\uvec_N\\ \vvec_1,\dots,\vvec_N\\ (\tau_{\vvec_1}<
 \tau_{\uvec_2})}}&\bigg(\mP_p\Big(\{\ovec\db\vvec_1\}\cap\{\ovec\to\uvec_2\}
 \Big)(q_p*\varphi_p)(\vvec_1;\vvec_2)\nn\\
&+\mP_p\Big(\{\ovec\db\vvec_1\}\cap\{\ovec\to\vvec_2\}\Big)(q_p*\varphi_p)
 (\vvec_1;\uvec_2)\bigg)\varphi_p(\uvec_2;\vvec_2)\nn\\
&\hspace{-1pc}\times\prod_{i=2}^{N-1}\Xi_p(\uvec_i,\vvec_i;\uvec_{i+1},
 \vvec_{i+1})\,\Xi_p(\uvec_N,\vvec_N;\xvec,\xvec).
\end{align}
By \refeq{dbellsup} and the BK inequality and using the Markov property and
\refeq{varphi-trivineq} under the restriction $\tau_{\vvec_1}<\tau_{\uvec_2}$,
we obtain \refeq{piN-BKappl}. \QED

\begin{lemma}\label{lem:PiN-BKappl}
For $N\ge1$,
\begin{align}\lbeq{PiN-BKappl}
\Pi_p^{\sss(N)}(\xvec)&\equiv\sum_{\vec b_N,b}\sum_{j=1}^N\mP_p\Big(\tilde
 E_{\vec b_N}^{\sss(N)}(\xvec)\cap\big\{b=b_j\text{ or }b\in\piv(\tb_j,\bb_{j
 +1})\big\}\Big)\nn\\
&\le\sum_{\substack{\uvec_1,\dots,\uvec_{N+1}\\ \vvec_1,\dots,\vvec_{N
 +1}\\ (\uvec_{N+1}=\vvec_{N+1}=\xvec)}}\varphi_p(\uvec_1)\,\varphi_p(\uvec_1;
 \vvec_1)\,\varphi_p(\vvec_1)\sum_{j=1}^N\prod_{i\ne j}\Xi_p(\uvec_i,\vvec_i;
 \uvec_{i+1},\vvec_{i+1})\\
&\hspace{7pc}\times\Big(\Xi_p(\uvec_j,\vvec_j;\uvec_{j+1},\vvec_{j+1})+\Theta_p
 (\uvec_j,\vvec_j;\uvec_{j+1},\vvec_{j+1})+\Theta'_p(\uvec_j,\vvec_j;\uvec_{j+
 1},\vvec_{j+1})\Big),\nn
\end{align}
where the empty product $\prod_{i\ne j}\Xi_p(\uvec_i,\vvec_i;\uvec_{i+1},
\vvec_{i+1})$ for the case of $N=1$ is 1 by convention, and
\begin{align}
&\Theta_p(\uvec,\vvec;\uvec',\vvec')=\big(\thetap_p(\uvec,\vvec;\uvec',
 \vvec')+\thetac_p(\uvec,\vvec;\uvec',\vvec')\big)\varphi_p(\uvec';\vvec')
 /2^{\delta_{\uvec',\vvec'}},\\
&\begin{cases}
\thetap_p(\uvec,\vvec;\uvec',\vvec')=(q_p*\varphi_p)(\uvec;\uvec')\,(q_p*
 \varphi_p*q_p*\varphi_p)(\vvec;\vvec'),\\
\thetac_p(\uvec,\vvec;\uvec',\vvec')=(q_p*\varphi_p)(\uvec;\vvec')\,(q_p*
 \varphi_p*q_p*\varphi_p)(\vvec;\uvec'),
\end{cases}\\
&\Theta'_p(\uvec,\vvec;\uvec',\vvec')=(q_p*\varphi_p)(\uvec;\vvec')\,(q_p*
 \varphi_p)(\vvec;\uvec')\,(\varphi_p*q_p*\varphi_p)(\uvec';\vvec').
 \lbeq{Xipp}
\end{align}
\end{lemma}

\Proof{Proof} Since
\begin{align}\lbeq{PiN-fixedj}
\sum_{\vec b_N,b}\mP_p\Big(\tilde E_{\vec b_N}^{\sss(N)}(\xvec)\cap\big\{b=
 b_j\text{ or }b\in\piv(\tb_j,\bb_{j+1})\big\}\Big)=\pi_p^{\sss(N)}(\xvec)+
 \sum_{\vec b_N,b}\mP_p\Big(\tilde E_{\vec b_N}^{\sss(N)}(\xvec)\cap\big\{b
 \in\piv(\tb_j,\bb_{j+1})\big\}\Big),
\end{align}
it suffices to investigate the sum on the right-hand side.  To do so, we use
the following relations that are similar to \refeq{Esup} and \refeq{Eellsup}:
\begin{align}
E(b',\xvec;\tilde\cC^{b'}(\yvec))\cap\big\{b\in\piv(\tbp,\xvec)\big\}&\subset
 \{\yvec\to\xvec\}\circ\{b'\to b\to\xvec\},\lbeq{Epivsup}\\
E(b',\vvec;\tilde\cC^{b'}(\yvec))\cap\big\{b\in\piv(\tbp,\vvec)\big\}\cap\{\tbp
 \to\xvec\}&\subset\bigcup_{\uvec:\tau_{\uvec}>\tau_{\bb'}}\Big\{\big\{\{\yvec
 \to\uvec\to\vvec\}\circ\{b'\to b\to\vvec\}\circ\{\uvec\to\xvec\}\big\}\nn\\
&\hspace{3.7pc}\cup\big\{\{\yvec\to\vvec\}\circ\{b'\to b\to\uvec\to\vvec\}\circ
 \{\uvec\to\xvec\}\big\}\Big\}\nn\\
&\hspace{3.7pc}\cup\big\{\{\yvec\to\vvec\}\circ\{b'\to\uvec\to b\to\vvec\}\circ
 \{\uvec\to\xvec\}\big\}\Big\}.\lbeq{Eellpivsup}
\end{align}

First we let $j=N$.  By \refeq{Epivsup} and using the BK inequality and the
Markov property, we obtain
\begin{align}\lbeq{PiN-1stbd}
\sum_{b_N,b}\mP_p\Big(\tilde E_{\vec b_N}^{\sss(N)}(\xvec)\cap\big\{b\in\piv
 (\tb_N,\xvec)\big\}\Big)\le\sum_{\vvec_N}\mP_p\Big(\tilde E_{\vec b_{N-1}}
 ^{\sss(N-1)}(\vvec_N)\cap\{\tb_{N-1}\to\xvec\}\Big)(q_p*\varphi_p*q_p*
 \varphi_p)(\vvec_N;\xvec),
\end{align}
which is equivalent to \refeq{piN-1stbd}, except for the last term
$(q_p*\varphi_p*q_p*\varphi_p)(\vvec_N;\xvec)$.  Therefore, by following the
same line as in \refeq{piN-2ndbd}--\refeq{piN-5thbd}, we obtain
\begin{align}\lbeq{PiN-2ndbd}
\sum_{\vec b_N,b}\mP_p\Big(\tilde E_{\vec b_N}^{\sss(N)}(\xvec)\cap\big\{b
 \in\piv(\tb_N,\xvec)\big\}\Big)\le\sum_{\substack{\uvec_1,\dots,\uvec_N\\
 \vvec_1,\dots,\vvec_N}}\varphi_p(\uvec_1)\,\varphi_p(\uvec_1;\vvec_1)\,
 \varphi_p(\vvec_1)&\prod_{i=1}^{N-1}\Xi_p(\uvec_i,\vvec_i;\uvec_{i+1},
 \vvec_{i+1})\nn\\
&\times\Theta_p(\uvec_N,\vvec_N;\xvec,\xvec).
\end{align}
Applying \refeq{piN-BKappl} to $\pi_p^{\sss(N)}(\xvec)$ in \refeq{PiN-fixedj}
and using $\Theta'_p(\uvec,\vvec;\uvec',\vvec')=0$ for $\uvec'=\vvec'$ (since
$(\varphi_p*q_p*\varphi_p)(\uvec';\vvec')$ in \refeq{Xipp} is zero if
$\uvec'=\vvec'$), we obtain the term for $j=N$ in \refeq{PiN-BKappl}.

Next we let $j<N$.  Following the same line as in
\refeq{piN-1stbd}--\refeq{piN-4thbd}, we obtain
\begin{align}\lbeq{PiN-3rdbd}
&\sum_{\vec b_N,b}\mP_p\Big(\tilde E_{\vec b_N}^{\sss(N)}(\xvec)\cap\big\{b\in
 \piv(\tb_j,\bb_{j+1})\big\}\Big)\nn\\
&\le\sum_{\substack{\uvec_{j+2},\dots,\uvec_N\\ \vvec_{j+1},\dots,\vvec_N\\
 (\tau_{\vvec_{j+1}}<\tau_{\uvec_{j+2}})}}\sum_{\vec b_j,b}\bigg(\mP_p\Big(
 \tilde E_{\vec b_j}^{\sss(j)}(\vvec_{j+1})\cap\big\{b\in\piv(\tb_j,\vvec_{
 j+1})\big\}\cap\{\tb_j\to\uvec_{j+2}\}\Big)(q_p*\varphi_p)(\vvec_{j+1};
 \vvec_{j+2})\nn\\
&\hspace{7.5pc}+\mP_p\Big(\tilde E_{\vec b_j}^{\sss(j)}(\vvec_{j+1})\cap\big\{b
 \in\piv(\tb_j,\vvec_{j+1})\big\}\cap\{\tb_j\to\vvec_{j+2}\}\Big)(q_p*
 \varphi_p)(\vvec_{j+1};\uvec_{j+2})\bigg)\nn\\
&\hspace{6.5pc}\times\varphi_p(\uvec_{j+2};\vvec_{j+2})\prod_{i=j+2}^{N-1}\Xi_p
 (\uvec_i,\vvec_i;\uvec_{i+1},\vvec_{i+1})\,\Xi_p(\uvec_N,\vvec_N;\xvec,\xvec).
\end{align}
Then, by \refeq{Eellpivsup} and using the BK inequality and the Markov
property,
\begin{align}\lbeq{PiN-4thbd}
&\sum_{b_j,b}\mP_p\Big(\tilde E_{\vec b_j}^{\sss(j)}(\vvec_{j+1})\cap\big\{b
 \in\piv(\tb_j,\vvec_{j+1})\big\}\cap\{\tb_j\to\uvec_{j+2}\}\Big)\nn\\
&\le\sum_{\substack{\vvec_j,\uvec_{j+1}\\ (\tau_{\vvec_j}<\tau_{\uvec_{j+1}})}}
 \Bigg(\bigg(\mP_p\Big(\tilde E_{\vec b_{j-1}}^{\sss(j-1)}(\vvec_j)\cap\{
 \tb_{j-1}\to\uvec_{j+1}\}\Big)(q_p*\varphi_p*q_p*\varphi_p)(\vvec_j;
 \vvec_{j+1})\nn\\
&\hspace{6pc}+\mP_p\Big(\tilde E_{\vec b_{j-1}}^{\sss(j-1)}(\vvec_j)\cap\{
 \tb_{j-1}\to\vvec_{j+1}\}\Big)(q_p*\varphi_p*q_p*\varphi_p)(\vvec_j;
 \uvec_{j+1})\bigg)\varphi_p(\uvec_{j+1};\vvec_{j+1})\\
&\hspace{3pc}+\mP_p\Big(\tilde E_{\vec b_{j-1}}^{\sss(j-1)}(\vvec_j)\cap\{
 \tb_{j-1}\to\vvec_{j+1}\}\Big)(q_p*\varphi_p)(\vvec_j;\uvec_{j+1})\,
 (\varphi_p*q_p*\varphi_p)(\uvec_{j+1};\vvec_{j+1})\Bigg)\varphi_p(\uvec_{j
 +1};\uvec_{j+2}),\nn
\end{align}
where the last term can be replaced by
$(q_p*\varphi_p)(\uvec_{j+1};\uvec_{j+2})$, because
$\tau_{\uvec_{j+1}}\le\tau_{\vvec_{j+1}}<\tau_{\uvec_{j+2}}$ (due to the
restriction in \refeq{PiN-3rdbd} and the factors
$\varphi_p(\uvec_{j+1};\vvec_{j+1})$ and
$(\varphi_p*q_p*\varphi_p)(\uvec_{j+1};\vvec_{j+1})$ in \refeq{PiN-4thbd}).
Using \refeq{PiN-4thbd} as well as that with $\uvec_{j+2}$ replaced by
$\vvec_{j+2}$, we obtain
\begin{align}\lbeq{PiN-5thbd}
\refeq{PiN-3rdbd}&\le\sum_{\substack{\uvec_{j+1},\dots,\uvec_N\\ \vvec_j,
 \dots,\vvec_N\\ (\tau_{\vvec_j}<\tau_{\uvec_{j+1}})}}\sum_{\vec b_{j-1}}
 \Bigg(\bigg(\mP_p\Big(\tilde E_{\vec b_{j-1}}^{\sss(j-1)}(\vvec_j)\cap\{
 \tb_{j-1}\to\uvec_{j+1}\}\Big)(q_p*\varphi_p*q_p*\varphi_p)(\vvec_j;
 \vvec_{j+1})\nn\\
&\hspace{7.2pc}+\mP_p\Big(\tilde E_{\vec b_{j-1}}^{\sss(j-1)}(\vvec_j)\cap\{
 \tb_{j-1}\to\vvec_{j+1}\}\Big)\,(q_p*\varphi_p*q_p*\varphi_p)(\vvec_j;
 \uvec_{j+1})\bigg)\varphi_p(\uvec_{j+1};\vvec_{j+1})\nn\\
&\hspace{6.2pc}+\mP_p\Big(\tilde E_{\vec b_{j-1}}^{\sss(j-1)}(\vvec_j)\cap\{
 \tb_{j-1}\to\vvec_{j+1}\}\Big)\,(q_p*\varphi_p)(\vvec_j;\uvec_{j+1})\,
 (\varphi_p*q_p*\varphi_p)(\uvec_{j+1};\vvec_{j+1})\Bigg)\nn\\
&\hspace{5.2pc}\times\prod_{i=j+1}^{N-1}\Xi_p(\uvec_i,\vvec_i;\uvec_{i
 +1},\vvec_{i+1})\,\Xi_p(\uvec_N,\vvec_N;\xvec,\xvec).
\end{align}
Repeatedly using \refeq{piN-2ndbd}--\refeq{varphi-trivineq} with different
variables, we finally arrive at
\begin{align}\lbeq{PiN-6thbd}
\refeq{PiN-5thbd}&\le\sum_{\substack{\uvec_j,\dots,\uvec_{N+1}\\ \vvec
 _{j-1},\dots,\vvec_{N+1}\\ (\uvec_{N+1}=\vvec_{N+1}=\xvec)\\ (\tau_{\vvec_{j
 -1}}<\tau_{\uvec_j})}}\sum_{\vec b_{j-2}}\bigg(\mP_p\Big(\tilde E_{\vec b_{j
 -2}}^{\sss(j-2)}(\vvec_{j-1})\cap\{\tb_{j-2}\to\uvec_j\}\Big)(q_p*\varphi_p*
 q_p*\varphi_p)(\vvec_{j-1};\vvec_j)\nn\\
&\hspace{8.3pc}+\mP_p\Big(\tilde E_{\vec b_{j-2}}^{\sss(j-2)}(\vvec_{j-1})\cap
 \{\tb_{j-2}\to\vvec_j\}\Big)\,(q_p*\varphi_p*q_p*\varphi_p)(\vvec_{j-1};
 \uvec_j)\bigg)\varphi_p(\uvec_j;\vvec_j)\nn\\
&\hspace{7.3pc}\times\Big(\Theta_p(\uvec_j,\vvec_j;\uvec_{j+1},\vvec_{j+1})+
 \Theta'_p(\uvec_j,\vvec_j;\uvec_{j+1},\vvec_{j+1})\Big)\prod_{i=j+1}^N\Xi_p
 (\uvec_i,\vvec_i;\uvec_{i+1},\vvec_{i+1})\nn\\
&~~\vdots\nn\\
&\le\sum_{\substack{\uvec_1,\dots,\uvec_{N+1}\\ \vvec_1,\dots,\vvec_{N
 +1}\\ (\uvec_{N+1}=\vvec_{N+1}=\xvec)}}\varphi_p(\uvec_1)\,\varphi_p(\vvec_1-
 \uvec_1)\,\varphi_p(\vvec_1)\prod_{i\ne j}\Xi_p(\uvec_i,\vvec_i;\uvec_{i+1},
 \vvec_{i+1})\nn\\
&\hspace{7pc}\times\Big(\Theta_p(\uvec_j,\vvec_j;\uvec_{j+1},\vvec_{j+1})+
 \Theta'_p(\uvec_j,\vvec_j;\uvec_{j+1},\vvec_{j+1})\Big).
\end{align}
Combining this with the bound \refeq{piN-BKappl} on $\pi_p^{\sss(N)}(\xvec)$ in
\refeq{PiN-fixedj}, we obtain the term for $j<N$ in \refeq{PiN-BKappl}.

The proof of \refeq{PiN-BKappl} is completed by summing the above bounds over
$j=1,\dots,N$. \QED

\section{Proof of \cite[Proposition~3.3]{cs?}}
In this section, we prove \cite[Proposition~3.3]{cs?} using
Lemmas~\ref{lem:piN-BKappl}--\ref{lem:PiN-BKappl} and assuming
translation-invariance.

Let $\varphi_p^{\sss(m)}(\vvec;\xvec)=\varphi_p(\vvec;\xvec)m^{\tau_{\xvec}-
\tau_{\vvec}}$. Recall that the weighted bubble $W_p^{\sss(m)}(k)$ and the
triangles $T_p^{\sss(m)}$ and $\tilde T_p$ are defined as
\begin{align}
W_p^{\sss(m)}(k)&=\sup_{\xvec\in\mZ^{d+1}}\sum_{\vvec}\big(1-\cos(k\cdot
 \sigma_{\vvec})\big)\times\begin{cases}
 (q_p*\varphi_p)(\vvec)\,(mq_p*\varphi_p^{\sss(m)})(\xvec;\vvec)&\qquad(m<1),\\
 (mq_p*\varphi_p^{\sss(m)})(\vvec)\,(q_p*\varphi_p)(\xvec;\vvec)&\qquad(m\ge1),
 \end{cases}\lbeq{W-def}\\
T_p^{\sss(m)}&=\sup_{\xvec\in\mZ^{d+1}}\sum_{\vvec}(q_p*\varphi_p*\varphi_p)
 (\vvec)\,(mq_p*\varphi_p^{\sss(m)})(\xvec;\vvec),\lbeq{T-def}\\
\tilde T_p&=\sup_{\xvec\in\mZ^{d+1}}\sum_{\vvec}(q_p*\varphi_p*q_p*\varphi_p)
 (\vvec)\,(q_p*\varphi_p)(\xvec;\vvec),\lbeq{tildeT-def}
\end{align}
and that the square $S_p^{\sss(m)}$ and the H-shaped diagram $H_p$ are defined
as (cf., \cite[Figure~2]{cs?})
\begin{align}
S_p^{\sss(m)}&=\sup_{\xvec\in\mZ^{d+1}}\sum_{\vvec}(q_p*\varphi_p*\varphi_p*
 \varphi_p)(\vvec)\,(mq_p*\varphi_p^{\sss(m)})(\xvec;\vvec),\lbeq{S-def}\\
H_p&=\sup_{\xvec,\yvec\in\mZ^{d+1}}\sum_{\uvec,\vvec,\wvec}(q_p*\varphi_p)
 (\uvec)\,(\varphi_p*q_p*\varphi_p)(\uvec;\vvec)\,(q_p*\varphi_p)(\xvec;
 \vvec)\,(q_p*\varphi_p)(\uvec;\wvec)\,(q_p*\varphi_p)(\vvec;\yvec+\wvec).
 \lbeq{H-def}
\end{align}
\cite[Proposition~3.3]{cs?} is an immediate consequence of the following lemma:

\begin{lemma}\label{lem:diagbd}
\begin{enumerate}[(i)]
\item
For $N\geq0$ and $\ell=0,1,2$,
\begin{align}
&\sum_{\xvec\in\Zd\times\mN}\tau_{\xvec}^\ell\,\pi_p^{\sss(N)}(\xvec)m^{\tau_{
 \xvec}}\le(N+1)^\ell(1+2T_p^{\sss(m)})(2T_p^{\sss(m)})^{(N-1)\vee0}\times
 \begin{cases}
 T_p^{\sss(m)}&(\ell\le1),\\
 S_p^{\sss(m)}&(\ell=2),
 \end{cases}\lbeq{piNn-diagbd}\\
&\sum_{\xvec\in\Zd\times\Zp}\big(1-\cos(k\cdot\sigma_{\xvec})\big)\pi_p^{\sss
 (N)}(\xvec)m^{\tau_{\xvec}}\le3(N+1)^2(1+2T_p^{\sss(m)})(2T_p^{\sss(m)})^{(N
 -1)\vee0}W_p^{\sss(m)}(k).\lbeq{piNcos-diagbd}
\end{align}
\item For $N\geq1$,
\begin{align}\lbeq{PiN-diagbd}
\sum_{\xvec\in\Zd\times\Zp}\Pi_p^{\sss(N)}(\xvec)\le N(1+2T_p^{\sss(1)})\Big(
 (T_p^{\sss(1)}+\tilde T_p)(2T_p^{\sss(1)})^{N-1}+H_p(2T_p^{\sss(1)})^{(N-2)
 \vee0}\Big).
\end{align}
\end{enumerate}
\end{lemma}

\Proof{Proof of Lemma~\ref{lem:diagbd}(i)} First we prove
\underline{\refeq{piNn-diagbd}--\refeq{piNcos-diagbd} for $N=0$}.  By
\refeq{pi0-BKappl}, we readily obtain
\begin{align}\lbeq{pi0n-diagbd}
\sum_{\xvec\in\Zd\times\mN}\tau_{\xvec}^\ell\,\pi_p^{\sss(0)}(\xvec)m^{\tau_{
 \xvec}}\le\sum_{\xvec}\tau_{\xvec}^\ell\,(q_p*\varphi_p)(\xvec)\,(mq_p*
 \varphi_p^{\sss(m)})(\xvec)\le
 \begin{cases}
 T_p^{\sss(m)}&(\ell\le1),\\
 S_p^{\sss(m)}&(\ell=2),
 \end{cases}
\end{align}
where we have used (cf., \cite[(5.17)]{S1})
\begin{align}
&\tau_{\xvec}\,(q_p*\varphi_p)(\xvec)=\sum_{t=1}^{\tau_{\xvec}}\,(q_p*
 \varphi_p)(\xvec)\le\sum_{t=1}^{\tau_{\xvec}}\sum_{\vvec:\tau_{\vvec}=t}(q_p
 *\varphi_p)(\vvec)\,\varphi_p(\vvec;\xvec)=(q_p*\varphi_p*\varphi_p)(\xvec),
 \lbeq{Markov}\\
&\tau_{\xvec}^2\,(q_p*\varphi_p)(\xvec)\stackrel{\because\refeq{Markov}}\le\tau
 _{\xvec}\,(q_p*\varphi_p*\varphi_p)(\xvec)\stackrel{\because\refeq{Markov}}\le
 (q_p*\varphi_p*\varphi_p*\varphi_p)(\xvec).\lbeq{Markov2}
\end{align}
We note that we have multiplied one of the two diagram lines (i.e.,
$(q_p*\varphi_p)(\xvec)$) by $\tau_{\xvec}^\ell$ and the other by
$m^{\tau_{\xvec}}$.  If we multiply either $(q_p*\varphi_p)(\xvec)$ or
$(mq_p*\varphi_p^{\sss(m)})(\xvec)$ (depending on whether $m<1$ or $m\ge1$) by
$1-\cos(k\cdot\sigma_{\xvec})$ instead of $\tau_{\xvec}^\ell$, we obtain
\begin{align}\lbeq{pi0cos-diagbd}
\sum_{\xvec}\big(1-\cos(k\cdot\sigma_{\xvec})\big)\,\pi_p^{\sss(0)}(\xvec)
 m^{\tau_{\xvec}}\le W_p^{\sss(m)}(k),
\end{align}
as required.

Next we prove \underline{\refeq{piNn-diagbd} for $N\ge1$ and $\ell=0$}.  We
note that, as in the $N=0$ case above, there are two ``external'' diagram lines
from $\ovec$ to $\xvec$ in each of the $2^{N-1}$ bounding diagrams in
\refeq{piN-BKappl}.  Each line looks like
\begin{align}
\varphi_p(\yvec_1)\prod_{i=1}^{N-1}(q_p*\varphi_p)(\yvec_i;\yvec_{i+1})\,
 (q_p*\varphi_p)(\yvec_N;\xvec),
\end{align}
where each $\yvec_i$ is either $\uvec_i$ or $\vvec_i$ in \refeq{piN-BKappl};
denote the line with $\yvec_1=\vvec_1$ by
$\omega_{\vvec_1}\equiv(\ovec,\vvec_1,\dots,\xvec)$ and the other by
$\omega_{\uvec_1}\equiv(\ovec,\uvec_1,\dots,\xvec)$.  Multiplying
$\omega_{\vvec_1}$ by $m^{\tau_{\xvec}}$ and using
\begin{align}\lbeq{xibd}
\left.\begin{array}{c}
\dpst2\sum_{\xvec}\Xi_p(\ovec,\yvec;\xvec,\xvec)m^{\tau_{\xvec}}\\[1pc]
\dpst\sum_{\uvec,\vvec}\big(\xip_p(\ovec,\yvec;\uvec,\vvec)m^{\tau_{\uvec}}
 +\xic_p(\ovec,\yvec;\uvec,\vvec)m^{\tau_{\vvec}}\big)\varphi_p(\uvec;\vvec)
 \\[1.2pc]
\dpst\sum_{\uvec,\vvec}\big(\xip_p(\yvec,\ovec;\uvec,\vvec)m^{\tau_{\vvec}}
 +\xic_p(\yvec,\ovec;\uvec,\vvec)m^{\tau_{\uvec}}\big)\varphi_p(\uvec;\vvec)
\end{array}\!\!\right\}\le2T_p^{\sss(m)}\qquad(\yvec\in\mZ^{d+1}),
\end{align}
and
\begin{align}\lbeq{T0bd}
\sum_{\uvec,\vvec}\varphi_p(\uvec)\,\varphi_p(\uvec;\vvec)\,\varphi_p(\vvec)
 m^{\tau_{\vvec}}\le1+\sum_{\vvec\ne\ovec}(\varphi_p*\varphi_p)(\vvec)\,(mq_p*
 \varphi_p^{\sss(m)})(\vvec)\le1+2T_p^{\sss(m)},
\end{align}
we obtain
\begin{align}\lbeq{ell=0-bd}
\sum_{\xvec}\pi_p^{\sss(N)}(\xvec)m^{\tau_{\xvec}}\le(1+2T_p^{\sss(m)})
 (2T_p^{\sss(m)})^{N-1}T_p^{\sss(m)}\qquad(N\ge1),
\end{align}
as required.

Before proceeding the proof, we define (cf., \refeq{xi-def})
\begin{align}\lbeq{tildexi-def}
\tilde\Xi_p(\uvec,\vvec;\uvec',\vvec')=\varphi_p(\uvec;\vvec)\big(\xip_p(\uvec,
 \vvec;\uvec',\vvec')+\xic_p(\uvec,\vvec;\uvec',\vvec')\big)/2^{\delta_{\uvec',
 \vvec'}},
\end{align}
which satisfies similar bounds to \refeq{xibd}, due to translation-invariance.
We note that, by using \refeq{tildexi-def}, the bound in \refeq{piN-BKappl} can
be reorganized as
\begin{align}\lbeq{piN-BKappl-rewr1}
\sum_{\substack{\uvec_1,\dots,\uvec_{N+1}\\ \vvec_1,\dots,\vvec_{N+1}\\
 (\uvec_{N+1}=\vvec_{N+1}=\xvec)}}\varphi_p(\uvec_1)\,\varphi_p(\vvec_1)
 \prod_{i=1}^N\tilde\Xi_p(\uvec_i,\vvec_i;\uvec_{i+1},\vvec_{i+1}),
\end{align}
or, for $j=1,\dots,N$, as
\begin{align}\lbeq{piN-BKappl-rewr2}
&\sum_{\substack{\uvec_1,\dots,\uvec_{N+1}\\ \vvec_1,\dots,\vvec_{N+1}\\
 (\uvec_{N+1}=\vvec_{N+1}=\xvec)}}\varphi_p(\uvec_1)\,\varphi_p(\uvec_1;
 \vvec_1)\,\varphi_p(\vvec_1)\bigg(\prod_{i=1}^{j-1}\Xi_p(\uvec_i,\vvec_i;
 \uvec_{i+1},\vvec_{i+1})\bigg)\nn\\
&\qquad\times\Big(\xip_p(\uvec_j,\vvec_j;\uvec_{j+1},\vvec_{j+1})+\xic_p
 (\uvec_j,\vvec_j;\uvec_{j+1},\vvec_{j+1})\Big)\bigg(\prod_{i=j+1}^N\tilde
 \Xi_p(\uvec_i,\vvec_i;\uvec_{i+1},\vvec_{i+1})\bigg).
\end{align}

Now we prove \underline{\refeq{piNn-diagbd} for $N\ge1$ and $\ell=1,2$}.  To do
so, we multiply $\omega_{\vvec_1}$ by $m^{\tau_{\xvec}}$ as before, and
multiply $\omega_{\uvec_1}=(\ovec,\omega_{\uvec_1}^{\sss
(1)},\dots,\omega_{\uvec_1}^{\sss(N+1)})$, where $\omega_{\uvec_1}^{\sss(1)}
=\uvec_1$, $\omega_{\uvec_1}^{\sss(i)}\in\{\uvec_i,\vvec_i\}$ for $i=2,\dots,N$
and $\omega_{\uvec_1}^{\sss(N+1)}=\xvec$, by $\tau_{\xvec}^\ell$, using the
decomposition
\begin{align}\lbeq{telescope-n}
\tau_{\xvec}=\tau_{\uvec_1}+\sum_{j=1}^N(\tau_{\omega_{\uvec_1}^{\sss(j+1)}}
 -\tau_{\omega_{\uvec_1}^{\sss(j)}}).
\end{align}
Consider, e.g., the bounding diagram with $\omega_{\uvec_1}^{\sss(i)}=\uvec_i$
for all $i=2,\dots,N$; we denote this diagram by $U(\xvec)$ for convenience.
Then, by \refeq{piN-BKappl-rewr1} and \refeq{Markov}--\refeq{Markov2}, the
contribution from $\tau_{\uvec_1}^\ell$ is bounded as
\begin{align}\lbeq{ell=1-bd1}
\sum_{\substack{\uvec_1,\dots,\uvec_{N+1}\\ \vvec_1,\dots,\vvec_{N+1}\\
 (\uvec_{N+1}=\vvec_{N+1})}}\!\!\!\tau_{\uvec_1}^\ell\varphi_p(\uvec_1)\,\varphi_p
 ^{\sss(m)}(\vvec_1)\prod_{i=1}^N\varphi_p(\uvec_i,\vvec_i)\,\xip_p(\uvec_i,
 \vvec_i;\uvec_{i+1},\vvec_{i+1})m^{\tau_{\vvec_{i+1}}-\tau_{\vvec_i}}\le
 (T_p^{\sss(m)})^N\times
 \begin{cases}
 T_p^{\sss(m)}&(\ell=1),\\
 S_p^{\sss(m)}&(\ell=2),
 \end{cases}
\end{align}
and, by \refeq{piN-BKappl-rewr2} and \refeq{Markov}--\refeq{Markov2} and using
\refeq{T0bd}, the contribution from each
$(\tau_{\uvec_{j+1}}-\tau_{\uvec_j})^\ell$ is bounded as
\begin{align}\lbeq{ell=1-bd2}
&\sum_{\substack{\uvec_1,\dots,\uvec_{N+1}\\ \vvec_1,\dots,\vvec_{N+1}\\
 (\uvec_{N+1}=\vvec_{N+1})}}\varphi_p(\uvec_1)\,\varphi_p(\uvec_1;\vvec_1)\,
 \varphi_p^{\sss(m)}(\vvec_1)\bigg(\prod_{i=1}^{j-1}\xip_p(\uvec_i,\vvec_i;
 \uvec_{i+1},\vvec_{i+1})m^{\tau_{\vvec_{i+1}}-\tau_{\vvec_i}}\varphi_p
 (\uvec_{i+1};\vvec_{i+1})\bigg)\nn\\
&\hspace{5pc}\times(\tau_{\uvec_{j+1}}-\tau_{\uvec_j})^\ell\,\xip_p(\uvec_j,
 \vvec_j;\uvec_{j+1},\vvec_{j+1})m^{\tau_{\vvec_{j+1}}-\tau_{\vvec_j}}\nn\\
&\hspace{5pc}\times\bigg(\prod_{i=j+1}^N\varphi_p(\uvec_i;\vvec_i)\,\xip_p
 (\uvec_i,\vvec_i;\uvec_{i+1},\vvec_{i+1})m^{\tau_{\vvec_{i+1}}-\tau_{\vvec_i}}
 \bigg)\nn\\
&\le(1+2T_p^{\sss(m)})(T_p^{\sss(m)})^{N-1}\times
 \begin{cases}
 T_p^{\sss(m)}&(\ell=1),\\
 S_p^{\sss(m)}&(\ell=2).
 \end{cases}
\end{align}
Therefore, for $\ell=1$,
\begin{align}\lbeq{ell=1-bd}
\sum_{\xvec}\tau_{\xvec}U(\xvec)m^{\tau_{\xvec}}\le N(1+2T_p^{\sss(m)})(T_p^{
 \sss(m)})^N+(T_p^{\sss(m)})^{N+1}\le(N+1)(1+2T_p^{\sss(m)})(T_p^{\sss(m)})^N.
\end{align}
The other $2^{N-1}-1$ bounding diagrams obey the same bound.  This completes
the proof of \refeq{piNn-diagbd} for $\ell\le1$.

The cross terms for $\ell=2$ can also be bounded similarly.  For example, the
contribution from $(\tau_{\uvec_{j'+1}}-\tau_{\uvec_{j'}})
(\tau_{\uvec_{j+1}}-\tau_{\uvec_j})$ with $j'<j$ is bounded, by using
\refeq{piN-BKappl-rewr2} (cf., \refeq{ell=1-bd2}), by
\begin{align}\lbeq{ell=2-bd1}
&(T_p^{\sss(m)})^{N-j'}\sum_{\substack{\uvec_1,\dots,\uvec_{j'+1}\\ \vvec_1,
 \dots,\vvec_{j'+1}}}\varphi_p(\uvec_1)\,\varphi_p(\uvec_1;\vvec_1)\,
 \varphi_p^{\sss(m)}(\vvec_1)\prod_{i=1}^{j'-1}\xip_p(\uvec_i,\vvec_i;
 \uvec_{i+1},\vvec_{i+1})m^{\tau_{\vvec_{i+1}}-\tau_{\vvec_i}}\varphi_p
 (\uvec_{i+1};\vvec_{i+1})\nn\\
&\hspace{8pc}\times(\tau_{\uvec_{j'+1}}-\tau_{\uvec_{j'}})\,\xip_p(\uvec_{j'},
 \vvec_{j'};\uvec_{j'+1},\vvec_{j'+1})m^{\tau_{\vvec_{j'+1}}-\tau_{\vvec_{j'}}}
 \varphi_p(\uvec_{j'+1};\vvec_{j'+1})\nn\\[5pt]
&\le(1+2T_p^{\sss(m)})(T_p^{\sss(m)})^{N-1}S_p^{\sss(m)}.
\end{align}
There are $N(N-1)-1$ more cross terms that obey the same bound.  There are $2N$
cross terms remaining, each of which is bounded by
$(T_p^{\sss(m)})^NS_p^{\sss(m)}$.  Therefore,
\begin{align}\lbeq{ell=2-bd2}
\sum_{\xvec}\tau_{\xvec}^2U(\xvec)m^{\tau_{\xvec}}&\le N^2(1+2T_p^{\sss(m)})
 (T_p^{\sss(m)})^{N-1}S_p^{\sss(m)}+(2N+1)(T_p^{\sss(m)})^NS_p^{\sss(m)}\nn\\
&\le(N+1)^2(1+2T_p^{\sss(m)})(T_p^{\sss(m)})^{N-1}S_p^{\sss(m)}.
\end{align}
The other $2^{N-1}-1$ bounding diagrams than $U(\xvec)$ obey the same bound.
This completes the proof of \refeq{piNn-diagbd}.

Finally we prove \underline{\refeq{piNcos-diagbd} for $N\ge1$}.  If $m<1$, then
we multiply $\omega_{\vvec_1}$ by $m^{\tau_{\xvec}}$ as before, and multiply
$\omega_{\uvec_1}$ by $1-\cos(k\cdot\sigma_{\xvec})$ and use the decomposition
(cf., \cite[(4.50)]{Slad06})
\begin{align}\lbeq{telescope-cos}
1-\cos(k\cdot\sigma_{\xvec})\le(2N+3)\bigg(1-\cos(k\cdot\sigma_{\uvec_1})
 +\sum_{j=1}^N\Big(1-\cos\big(k\cdot(\sigma_{\omega_{\uvec_1}^{\sss(j+1)}}
 -\sigma_{\omega_{\uvec_1}^{\sss(j)}})\big)\Big)\bigg).
\end{align}
For example, consider the bounding diagram $U(\xvec)$ again, where
$\omega_{\uvec_1}^{\sss(i)}=\uvec_i$ for all $i=2,\dots,N$.  Similarly to
\refeq{ell=1-bd1}--\refeq{ell=1-bd2}, we have
\begin{align}
\sum_{\substack{\uvec_1,\dots,\uvec_{N+1}\\ \vvec_1,\dots,\vvec_{N+1}\\
 (\uvec_{N+1}=\vvec_{N+1})}}&\big(1-\cos(k\cdot\sigma_{\uvec_1})\big)\varphi_p
 (\uvec_1)\,\varphi_p^{\sss(m)}(\vvec_1)\prod_{i=1}^N\varphi_p(\uvec_i,\vvec_i)
 \,\xip_p(\uvec_i,\vvec_i;\uvec_{i+1},\vvec_{i+1})m^{\tau_{\vvec_{i+1}}-\tau_{
 \vvec_i}}\nn\\
&\le W_p^{\sss(m)}(k)(T_p^{\sss(m)})^N,\lbeq{cos-bd1}
\end{align}
and
\begin{align}\lbeq{cos-bd2}
&\sum_{\substack{\uvec_1,\dots,\uvec_{N+1}\\ \vvec_1,\dots,\vvec_{N+1}\\
 (\uvec_{N+1}=\vvec_{N+1})}}\varphi_p(\uvec_1)\,\varphi_p(\uvec_1;\vvec_1)\,
 \varphi_p^{\sss(m)}(\vvec_1)\bigg(\prod_{i=1}^{j-1}\xip_p(\uvec_i,\vvec_i;
 \uvec_{i+1},\vvec_{i+1})m^{\tau_{\vvec_{i+1}}-\tau_{\vvec_i}}\varphi_p
 (\uvec_{i+1};\vvec_{i+1})\bigg)\nn\\
&\hspace{5pc}\times\Big(1-\cos\big(k\cdot(\sigma_{\uvec_{j+1}}-\sigma_{
 \uvec_j})\big)\Big)\xip_p(\uvec_j,\vvec_j;\uvec_{j+1},\vvec_{j+1})m^{\tau_{
 \vvec_{j+1}}-\tau_{\vvec_j}}\nn\\
&\hspace{5pc}\times\bigg(\prod_{i=j+1}^N\varphi_p(\uvec_i;\vvec_i)\,\xip_p
 (\uvec_i,\vvec_i;\uvec_{i+1},\vvec_{i+1})m^{\tau_{\vvec_{i+1}}-\tau_{\vvec_i}}
 \bigg)\nn\\
&\le(1+2T_p^{\sss(m)})(T_p^{\sss(m)})^{N-1}W_p^{\sss(m)}(k).
\end{align}
Therefore,
\begin{align}
\sum_{\xvec}(1-\cos(k\cdot\sigma_{\xvec}))U(\xvec)m^{\tau_{\xvec}}&\le(2N+3)
 \Big((T_p^{\sss(m)})^NW_p^{\sss(m)}(k)+N(1+2T_p^{\sss(m)})(T_p^{\sss(m)})^{N
 -1}W_p^{\sss(m)}(k)\Big)\nn\\
&\le3(N+1)^2(1+2T_p^{\sss(m)})(T_p^{\sss(m)})^{N-1}W_p^{\sss(m)}(k).
\end{align}
The other $2^{N-1}-1$ bounding diagrams than $U(\xvec)$ obey the same
bound.

If $m\ge1$, then we multiply $\omega_{\uvec_1}$ by
$(1-\cos(k\cdot\sigma_{\xvec}))m^{\tau_{\xvec}}$ and use the decomposition
\refeq{telescope-cos}.  The rest is the same.  This completes the proof of
\refeq{piNcos-diagbd} for $N\ge1$. \QED

\Proof{Proof of Lemma~\ref{lem:diagbd}(ii)} First we recall \refeq{PiN-BKappl}.
Since we have the bound \refeq{ell=0-bd} on the contribution from
$\Xi_p(\uvec_j,\vvec_j;\uvec_{j+1},\vvec_{j+1})$, it thus remains to
investigate the contributions from
$\Theta_p(\uvec_j,\vvec_j;\uvec_{j+1},\vvec_{j+1})$ and
$\Theta'_p(\uvec_j,\vvec_j;\uvec_{j+1},\vvec_{j+1})$.  However, since (cf.,
\refeq{piN-BKappl-rewr2})
\begin{align}\lbeq{pc1}
&\sum_{\substack{\uvec_1,\dots,\uvec_{N+1}\\ \vvec_1,\dots,\vvec_{N+1}\\
 (\uvec_{N+1}=\vvec_{N+1})}}\varphi_p(\uvec_1)\,\varphi_p(\uvec_1;
 \vvec_1)\,\varphi_p(\vvec_1)\bigg(\prod_{i=1}^{j-1}\Xi_p(\uvec_i,\vvec_i;
 \uvec_{i+1},\vvec_{i+1})\bigg)\bigg(\prod_{i=j+1}^N\tilde\Xi_p(\uvec_i,
 \vvec_i;\uvec_{i+1},\vvec_{i+1})\bigg)\nn\\
&\hspace{6pc}\times\Big(\thetap_p(\uvec_j,\vvec_j;\uvec_{j+1},\vvec_{j+1})
 +\thetac_p(\uvec_j,\vvec_j;\uvec_{j+1},\vvec_{j+1})\Big)\le(1+2T_p^{\sss
 (1)})(2T_p^{\sss(1)})^{N-1}\tilde T_p,
\end{align}
and, for $j<N$,
\begin{align}\lbeq{pc2}
&\sum_{\substack{\uvec_1,\dots,\uvec_{N+1}\\ \vvec_1,\dots,\vvec_{N+1}\\
 (\uvec_{N+1}=\vvec_{N+1})}}\varphi_p(\uvec_1)\,\varphi_p(\uvec_1;
 \vvec_1)\,\varphi_p(\vvec_1)\bigg(\prod_{i=1}^{j-1}\Xi_p(\uvec_i,\vvec_i;
 \uvec_{i+1},\vvec_{i+1})\bigg)\bigg(\prod_{i=j
 +2}^N\tilde\Xi_p(\uvec_i,\vvec_i;\uvec_{i+1},\vvec_{i+1})\bigg)\nn\\
&\hspace{6pc}\times\Theta'_p(\uvec_j,\vvec_j;\uvec_{j+1},\vvec_{j+1})\Big(
 \xip_p(\uvec_{j+1},\vvec_{j+1};\uvec_{j+2},\vvec_{j+2})+\xic_p(\uvec_{j+1},
 \vvec_{j+1};\uvec_{j+2},\vvec_{j+2})\Big)\nn\\
&\quad\le(1+2T_p^{\sss(1)})(2T_p^{\sss(1)})^{N-2}H_p,
\end{align}
we obtain
\begin{align}
\sum_{\xvec}\Pi_p^{\sss(N)}(\xvec)&\le(1+2T_p^{\sss(1)})\Big(N(T_p^{\sss(1)}
 +\tilde T_p)(2T_p^{\sss(1)})^{N-1}+(N-1)H_p(2T_p^{\sss(1)})^{N-2}\Big)\nn\\
&\le N(1+2T_p^{\sss(1)})\Big((T_p^{\sss(1)}+\tilde T_p)(2T_p^{\sss(1)})^{N-1}
 +H_p(2T_p^{\sss(1)})^{(N-2)\vee0}\Big).
\end{align}
This completes the proof of Lemma~\ref{lem:diagbd}(ii). \QED

\end{document}